\documentclass[amsmath,english,a4paper,graphicx,12pt]{article}
\usepackage{doc}
\usepackage{graphicx}
\usepackage{amsmath}
\usepackage{amsfonts}
\usepackage{amssymb}
\usepackage{babel}
\usepackage{latexsym}
\usepackage{mathrsfs}

\newcommand{\eps}{\mbox{$\varepsilon$}}

\newcommand{\C}{\mathbb{C}}
\newcommand{\R}{\mathbb{R}}
\newcommand{\li}{\left}
\newcommand{\ri}{\right}
\newcommand{\abs}[1]{\left|#1\right|}

\title{\bf A fresh approach to the Paley--Wiener theorem for Mellin transforms and the Mellin--Hardy spaces}
\author{Carlo Bardaro, \thanks{
Department of Mathematics and Computer Sciences, University of Perugia,
via Vanvitelli 1, I-06123 Perugia, Italy, e-mail: 
carlo.bardaro@unipg.it} \and
Paul L. Butzer, \thanks{Lehrstuhl A fuer Mathematik, RWTH Aachen, Templergraben 55, Aachen, D-52056, Germany, e-mail: butzer@rwth-aachen.de}\and
Ilaria Mantellini, \thanks{Department of Mathematics and Computer Sciences, University of Perugia,
via Vanvitelli 1, I-06123 Perugia, Italy, e-mail: 
mantell@dmi.unipg.it}\and Gerhard Schmeisser\thanks{Department Mathematik, FAU Erlangen-N\"{u}rnberg, Cauerstr. 11, 91058 Erlangen, Germany, 
email: schmeisser@mi.uni-erlangen.de}
}
\begin{document}
\maketitle
\noindent
{\bf Abstract.} Here we give a new approach to the Paley--Wiener theorem in a Mellin analysis setting which avoids the use of the Riemann surface of the logarithm
 and analytical branches and is based on  new concepts of {\it polar-analytic function} in the Mellin setting and Mellin--Bernstein spaces. A notion of Hardy spaces in the Mellin setting is also given along with applications to exponential sampling  formulas of optical physics.

\section{\bf Introduction}

The structure of the Paley--Wiener space of all continuous functions $f \in L^2(\mathbb{R})$ having compactly supported Fourier transform 
is precisely described by the classical Paley--Wiener theorem of Fourier analysis. This basic result characterizes the Paley--Wiener spaces 
by the Bernstein spaces comprising all functions $f \in L^2(\mathbb{R})$ which have an analytic extension to the whole complex plane and 
are of exponential type (the type being connected with the bandwidth of the Fourier transform); see, e.g., \cite{BO}, \cite{YO}, \cite{RU}, 
\cite{RU1}, \cite{HIG}. Moreover, it  has a wide range of applications, especially in sampling theory and related fields; see 
\cite{HIG}, \cite{RS}, \cite{SCH}. Numerous variants of this theorem have been proved by several authors \cite{AN}, \cite{AD}, \cite{TU}. 
In particular, analogues of the Paley--Wiener theorem were obtained for other integral transforms. Recently, in \cite{BBMS} we have proved 
a version for Mellin transforms, by introducing the notion of a Mellin--Bernstein space, comprising all  functions 
$f \in X^2_c:=\{f:\mathbb{R}^+ \rightarrow \mathbb{C}: (\cdot) (\cdot)^{c -1/2} \in L^2(\mathbb{R}^+)\}$  which have an analytic extension to 
the Riemann surface of the (complex) logarithm and satisfy some exponential-type condition. We gave two different approaches, one involving 
purely complex analysis arguments and the other one using ``real'' arguments based on the statement of a Mellin extension of the classical 
Bernstein inequality also proved in \cite{BBMS}. Later on, in \cite{BBMS2} we applied our Paley--Wiener theorem in the Mellin setting to the 
study of the so-called Mellin distance of functions belonging to certain functional spaces (like Lipschitz spaces, Mellin--Sobolev spaces, and so on) 
from the Mellin--Bernstein space. This leads to precise estimates of the approximation error in certain basic formulae valid in Mellin--Bernstein 
spaces such as the exponential sampling formula and the Mellin reproducing kernel formula; see, e.g., \cite{BP}, \cite{BJ3}, \cite{BBM1}, \cite{BBM3}. 
The results in \cite{BBMS2} extend the corresponding ones in Fourier analysis (see \cite{BSS}, \cite{BSS2}) to the Mellin frame.

This paper is concerned with equivalent formulations of the 
Paley--Wiener theorem in the Mellin setting, its content being fully 
different from our papers \cite{BBMS} and \cite{BBMS2}. It is a new and simpler 
approach as it avoids the use of an abstract Riemann surface and analytic 
branches. It employs the helicoidal surface as a model of the Riemann surface of the (complex) logarithm and a notion of  analytic functions on it.
These considerations lead us to a concept of ``polar analyticity'' which 
enables us to introduce in a simple way  Bernstein classes and  
Hardy spaces in a Mellin setting.  Moreover, the proofs of the main theorems are notably
different from those of \cite{BBMS}. Concerning the notion of polar analyticity of a function $f$ at a point 
$(r_0, \theta_0) \in \mathbb{H}:=\{(r,\theta) \in \mathbb{R^+}\times \mathbb{R}\},$ here discussed, we introduce it by means of the limit 
(see Sec.\ 3, Definition 3)
$$\lim_{(r,\theta) \rightarrow (r_0,\theta_0)} \frac{f(r, \theta) - f(r_0, \theta_0)}{re^{i\theta} - r_0e^{i\theta_0}}.$$
Equivalently, polar-analytic functions can be described by the well-known Cauchy--Riemann equations with respect to polar coordinates; 
see, e.g., \cite{NP}, \cite{CBV}. The classical Hardy spaces are treated in many books; see, e.g., \cite{RU}, \cite{RR}, \cite{FS}.

In Section 2, we give some basic notions concerning Mellin analysis and preliminary results, including the Paley--Wiener theorem in the Mellin 
setting stated in \cite{BBMS}. In  Subsection 2.1, we illustrate the equivalent formulation in terms of the helicoidal surface. In Sections 3 and 4, 
we study polar-analytic functions and Mellin--Bernstein classes. The Paley--Wiener theorem is stated in Section 5. Sections 6 and 7 are devoted 
to the Hardy-type spaces in the Mellin frame and the study of the Mellin distance of a function $f$ belonging to a Hardy space from the 
Mellin--Bernstein class. In the final section, we apply the results to estimates of the approximation error in the exponential sampling  formula.

\section{\bf Basic notions and preliminary results}

Let $C(\mathbb{R}^+)$ be the space of all continuous functions defined on $\mathbb{R}^+,$ and $C^{(r)}(\mathbb{R}^+)$ be the space of all functions 
in $C(\mathbb{R}^+)$ with a derivative of order $r$ in $C(\mathbb{R}^+).$ Analogously, by $C^\infty(\mathbb{R}^+)$ we denote the space of all 
infinitely differentiable functions.
By $L^1_{\rm{loc}}(\mathbb{R}^+)$, we denote the space of all measurable functions which are integrable on every bounded interval in $\mathbb{R}^+.$

For $1\leq p < +\infty,$ let $L^p(\mathbb{R}^+)$~ be the space of all Lebesgue measurable and $p$-integrable complex-valued functions defined on 
$\mathbb{R}^+$ endowed with the usual norm $\|f\|_p.$ Analogous notations hold for functions 
defined on $\mathbb{R}.$

For $p=1$ and $c \in \mathbb{R},$ let us consider the space (see \cite{BJ2})
$$X_c = \{ f: \mathbb{R}^+\rightarrow \mathbb{C}: f(\cdot) (\cdot)^{c-1}\in L^1(\mathbb{R}^+) \}$$
endowed with the norm
$$ \| f\|_{X_c} := \|f(\cdot) (\cdot)^{c-1} \|_1 = \int_0^{+\infty} |f(u)|u^{c-1} du.$$

More generally, let $X^p_c$  denote the space of all  functions $f: \mathbb{R}^+\rightarrow \mathbb{C}$ such that 
$f(\cdot) (\cdot)^{c-1/p}\in L^p (\mathbb{R}^+)$ with $1<p< \infty.$ In an equivalent form, $X^p_c$ is the space of all functions $f$ such that 
$(\cdot)^c f(\cdot) \in L^p_\mu(\mathbb{R}^+),$ where $ L^p_\mu(\mathbb{R}^+)$ denotes the Lebesgue space with respect to the (invariant) measure 
$\mu (A) = \int_A dt/t$ for  any measurable set $A \subset \mathbb{R}^+.$ Finally, for $p=\infty$, we define $X^\infty_c$ as the space comprising 
all measurable functions $f : \mathbb{R}^+\rightarrow \mathbb{C}$ such that $\|f\|_{X^\infty_c}:= \sup_{x>0}x^{c}|f(x)| < \infty.$ For $p=2$ see \cite{BJ4}.

 The Mellin translation operator $\tau_h^c$, for $h \in \mathbb{R}^+,~c \in \mathbb{R},$~$f: \mathbb{R}^+ \rightarrow \mathbb{C},$ is denoted by
$$(\tau_h^c f)(x) := h^c f(hx)~~(x\in \mathbb{R}^+).$$
Setting $\tau_h:= \tau^0_h,$ we have $(\tau_h^cf)(x) = h^c (\tau_hf)(x)$ and  $\|\tau_h^c f\|_{X_c} = \|f\|_{X_c}.$
\vskip0.3cm
In the Mellin frame, the natural concept of a  pointwise derivative of a function $f$ is given by the limit of the difference quotient involving 
the Mellin translation; thus if $f'$ exists,
$$\lim_{h \rightarrow 1}\frac{\tau_h^cf(x) - f(x)}{h-1} =  x f'(x) + cf(x).$$
This gives the motivation for the following definition (see \cite{BJ2}):
The pointwise Mellin differential operator $\Theta_c,$ or the pointwise Mellin derivative $\Theta_cf$ of a function 
$f: \mathbb{R}^+ \rightarrow \mathbb{C}$ and $c \in \mathbb{R},$ is defined by
\begin{eqnarray*}
\Theta_cf(x) := x f'(x) + c f(x)~~~(x \in \mathbb{R}^+)
\end{eqnarray*}
provided that $f'$ exists a.e.\ on $\mathbb{R}^+.$ The Mellin differential operator of order $r \in \mathbb{N}$ is defined recursively by
\begin{eqnarray*}
\Theta^1_c := \Theta_c ,\quad\quad \Theta^r_c := \Theta_c (\Theta_c^{r-1}).
\end{eqnarray*}
For convenience, set $\Theta^r:= \Theta^r_0$ for $c=0$ and $\Theta_c^0 := I$ with $I$ denoting the identity operator.
For instance, the first three Mellin derivatives are given by:
\begin{eqnarray*}
\Theta_cf(x) &=& xf'(x) + cf(x),\\
\Theta^2_cf(x) &=& x^2 f''(x) + (2c+1) xf'(x) + c^2f(x),\\
\Theta^3_cf(x) &= &x^3 f'''(x) + (3c+3)x^2f''(x) \\&& +(3c^2 + 3c +1)xf'(x) + c^3 f(x).
\end{eqnarray*}

The Mellin transform of a function $f\in X_c$ is the linear and bounded operator defined by (see, e.g., \cite{MA}, \cite{GPS},  \cite{BJ2})
$$ M_c[f](s) \equiv [f]^{\wedge}_{M_c} (s) := \int_0^{+\infty} u^{s-1} f(u) du~~~(s=c+ it, t\in \mathbb{R}).$$
The inverse Mellin transform $M^{-1}_c[g]$ of a function $g \in L^1(\{c\} \times i \mathbb{R}),$ is defined by:
\begin{eqnarray*}
M^{-1}_c[g](x)  :=
 \frac{x^{-c}}{2 \pi}\int_{-\infty}^{+\infty} g(c+it) x^{-it}dt ~~~(x \in \mathbb{R}^+),
\end{eqnarray*}
where  in general $L^p(\{c\} \times i \mathbb{R}),$ for $p \geq 1,$ will mean the space of all functions 
$g:\, c+i \mathbb{R} \rightarrow \mathbb{C}$ with 
$g(c +i\cdot) \in L^p(\mathbb{R}^+).$
\vskip0,3cm
We have the following preliminary results (see \cite{BJ2}, \cite{BBM1}):
\vskip0,3cm
\newtheorem{Lemma}{Lemma}
\begin{Lemma}[Inversion Theorem in $X_c$] \label{inversion}
If $f \in X_c$ is such that $M_c[f] \in L^1(\{c\} \times i \mathbb{R}),$ then 
$$M_c^{-1}[M_c[f]](x) = \frac{x^{-c}}{2 \pi} \int_{-\infty}^{+\infty} [f]^\wedge_{M_c}(c+it)x^{-it}dt = f(x) \quad \quad (\hbox{a.e.\ on } \mathbb{R}^+).$$
\end{Lemma}
\vskip0,3cm
\noindent
 The following lemma will enable us to work  in a practical Hilbert space setting.
\begin{Lemma}\label{inclusion}
 If $f \in X_c$  and $M_c[f] \in L^1(\{c\} \times i \mathbb{R}),$ then 
$f \in X^2_c.$ 
\end{Lemma}
\vskip0,4cm

More generally, for $1<p \leq 2,$  the Mellin transform $M_c^p$ of $f \in X^p_c$ is given by (see \cite{BJ4})
$$M_c^p[f](s) \equiv [f]^{\wedge}_{M_c^p} (s) = \mbox{l.i.m.}_{\rho \rightarrow +\infty}~\int_{1/\rho}^\rho f(u) u^{s-1}du,$$
for $s=c+it,$ 
in the sense that
$$\lim_{\rho \rightarrow \infty}\bigg\|M_c^p[f](c+it) - \int_{1/\rho}^\rho f(u) u^{s-1}du\bigg\|_{L^{p'}(\{c\}\times i \mathbb{R})} = 0,$$
where $p'$ is the conjugate exponent of $p$, that is, $1/p + 1/p'=1$.
In the following we are interested in the case $p=2.$ 

Analogously, we define the inverse Mellin transform of a function $g \in X^2_c$ by 
$$M^{2, -1}_c[g](x) = \mbox{l.i.m.}_{\rho \rightarrow +\infty}~\frac{1}{2 \pi}\int_{1/\rho}^\rho g(c+it) x^{-c-it}dt,$$
and for any $f \in X^2_c$, there holds
$$M^{2, -1}_c[M_c^2[f]](x) = f(x)\quad \quad (\hbox{a.e.\ on } \mathbb{R}^+);$$
see \cite{BJ4}.

For functions in $X_c \cap X^2_c,$ we have the following important ``consistency'' property of the Mellin transform (see \cite{BJ4}):
\begin{Lemma} \label{consistency}
If $f \in X_c \cap X^2_c$, then the Mellin transforms $M_c[f]$ and $M_c^2[f]$ coincide, i.e., $M_c[f](c+it) = M_c^2[f](c+it)$ for almost all $t \in \mathbb{R}.$
\end{Lemma}
Moreover, the following Mellin version of the Plancherel Theorem holds (see \cite[Lemma 2.6]{BJ4}):
\begin{Lemma}\label{Plancherel}
The operator $M_c^2$ from $X^2_c$ onto $L^2(\{c\}\times i \mathbb{R})$ is bounded and norm preserving, i.e., for $f \in X^2_c$, we have
$$\|f\|_{X^2_c}= \frac{1}{\sqrt{2\pi}}\|M_c^2[f]\|_{L^2(\{c\}\times i \mathbb{R})}.$$
\end{Lemma}


\subsection{The Mellin--Bernstein spaces via Riemann surfaces}

We begin by introducing the following spaces:
 \newtheorem{Definition}{Definition}
\begin{Definition}\label{def1}
Let $B^{2}_{c, T}$ denote the space of all functions  $f\in X_c^2\cap C(\mathbb{R}^+)$ 
such that $[f]^\wedge_{M_c^2}(c+it) = 0$ a.e. for $|t| > T.$ Analogously, by $B^1_{c,T}$ we denote the space of all functions $f\in X_c\cap C(\mathbb{R}^+)$ 
such that $[f]^\wedge_{M_c}(c+it) = 0$ for all $|t| > T.$
\end{Definition}
These spaces are called {\it Paley--Wiener} spaces.
By Lemma \ref{inclusion}, we have $B^1_{c,T} \subset X^2_c$, and 
it is easily seen that $B^1_{c,T} \subset B^2_{c, T}.$ 

In \cite{BBMS} we have seen that a Mellin bandlimited function cannot be extended to the whole complex plane as an entire function, but 
we have proved  that it has an analytic extension on the Riemann surface $S_{{\rm log}}$ of the logarithm. 

Consider a function $g$ analytic on $S_{{\rm log}}.$ Then $g$ can be split into branches $g_k~(k \in \mathbb{Z})$ that are analytic on the slit 
complex plane $\Omega:=\C\setminus \R^+$. Furthermore, they are connected on the slit by the following properties:
\begin{enumerate}
\item[(i)] For $x>0$ the limits
$$ g_k^+(x):= \lim_{\varepsilon \rightarrow 0+} g_k(x+i\varepsilon)
\quad \hbox{and}\quad
g_k^-(x):= \lim_{\varepsilon \rightarrow 0+} g_k(x-i\varepsilon)$$
 exist and
\begin{eqnarray}\label{slit}
g_k^-(x)= g_{k+1}^+(x), \quad g_0^+(x)= g(x)
\quad (x>0).
\end{eqnarray}
\item[(ii)] For $x>0$, let $U_x$ be an open disk in the right half-plane
with center at $x$. Then $\psi_k\,:\, U_x\rightarrow \mathbb{C}$ with
\begin{eqnarray*}
 \psi_k(z):= \left\{
\begin{array}{lll}
g_k(z) & \hbox{  for } & z\in U_x, \, \Im z <0,\\
g_k^-(z) & \hbox{ for } & z\in U_x \cap\mathbb{R},\\
g_{k+1}(z) & \hbox{ for } & z\in U_x, \,  \Im z>0
\end{array}
\right.
\end{eqnarray*}
is analytic.
\end{enumerate}
The Mellin--Bernstein space $\widetilde{B}^2_{c,T}$ comprises all functions $f\in X_c^2$
for which $g(x):=x^cf(x)$ has an analytic extension on $S_{{\rm log}}$ with branches $g_k$ on $\Omega$ satisfying the following (additional) conditions:
\begin{enumerate}
\item[(iii)] There exists a constant $C>0$ such that for all $k\in\mathbb{Z}$ and
$\theta\in [0,2\pi]$
\begin{eqnarray}\label{exp}
|g_k(re^{i\theta})|\,\leq \, C e^{T|2\pi k +\theta|} \qquad (r>0).
\end{eqnarray}
\item[(iv)]
For $\theta\in [0,2\pi]$, we have
\begin{eqnarray}\label{unif}
\lim_{r\rightarrow 0} g_k(re^{i\theta})\,=\, \lim_{r\rightarrow \infty}
g_k(re^{i\theta})\,=\,0
\end{eqnarray}
uniformly with respect to $\theta$.
\end{enumerate}
In (\ref{exp}) and (\ref{unif}) the value $g_k(re^{i\theta})$ has to be
defined as $g_k^+(r)$ when $\theta=0$ and as $g_k^-(r)$ when
$\theta=2\pi$.
\vskip0,3cm
Now the Paley--Wiener theorem for the Mellin transform can be
stated as follows (see \cite{BBMS}):
\newtheorem{Theorem}{Theorem}
\begin{Theorem}[Paley--Wiener] \label{pw0}
$\quad \widetilde{B}^2_{c,T} \,=\,
B^2_{c,T}$.
\end{Theorem}
 \noindent
A simple and convenient model for $S_{{\rm log}}$ is the helicoidal surface in $\mathbb{R}^3$ defined by
$$ \mathbb{E}:= \{ (x_1,x_2,x_3) \in \mathbb{R}^3 : x_1= r \cos \theta,~ x_2= r\sin \theta,~x_3=\theta,~ r>0, ~\theta\in \mathbb{R}\}.$$
The subset obtained by setting $\theta =0$ on the right-hand side can be interpreted as $\mathbb{R}^+.$
Just as $\mathbb{C}$ is an extension of  $\mathbb{R},$ we shall see that  $\mathbb{E}$ takes a corresponding role for $\mathbb{R}^+$ in Mellin analysis.
\begin{figure}[htbp] 
\begin{center}
\includegraphics[scale=0.80]{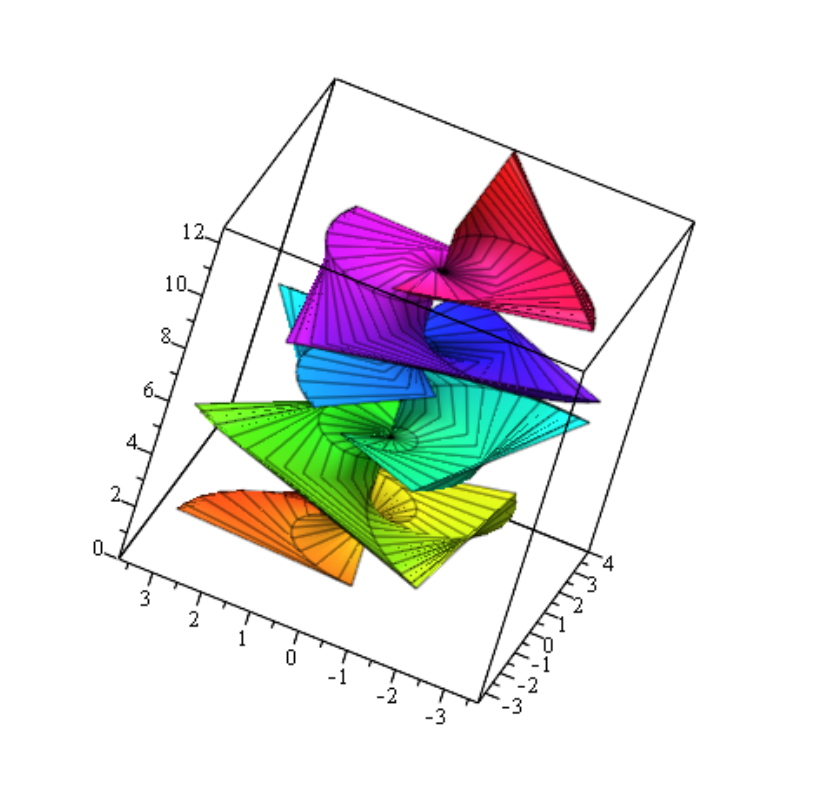}
\hskip0.2cm
\end{center}
\caption{\small The helicoidal surface as a model of the Riemann surface of the logarithm.}\label{helix}
\end{figure}
For $\alpha,~\beta\in \mathbb{R}$ with $\alpha<\beta,$ we consider the surface
$$\mathbb{E}_{\alpha,\beta}:= \left\{(x_1,x_2,x_3) \in \mathbb{R}^3 : x_1= r \cos \theta,~ x_2= r\sin \theta,~x_3=\theta,~ r>0, ~\theta\in {]}\alpha,\beta{[}\right\}$$
and call it a segment of $\mathbb{E}.$ The projection of $\mathbb{E}_{\alpha,\beta}$ into the $(x_1,x_2)$-plane may be interpreted as a sector in $\mathbb{C},$ given by
$$ S_{\alpha,\beta}:= \left\{ x_1 +ix_2 \in \mathbb{C} : x_1= r\cos\theta,~ x_2=r\sin \theta,~ r>0,~ \theta\in {]}\alpha,\beta{[} \right\}.$$
When $\beta -\alpha \in~{]}0, 2\pi],$ then this projection is a bijection. Indeed, given $z\in  S_{\alpha,\beta},$ there exists a unique $\theta\in~ {]}\alpha,\beta{[},$ denoted by $\theta:= \mbox{arg}_\alpha z,$ such that $z= |z| e^{i {\rm arg}_\alpha z}$ and then  $(\Re z, \Im z  ,\mbox{arg}_\alpha z)$ is the pre-image of $z$ on $\mathbb{E}_{\alpha,\beta}.$
Now, equivalently to the usual abstract approach (see \cite{BBMS}), analytic functions on $ \mathbb{E}$ may be introduced as follows.
\begin{Definition}\label{def3}
A function $f: \mathbb{E} \rightarrow \mathbb{C}$ is said to be {\rm analytic} if for every segment 
$\mathbb{E}_{\alpha,\beta}$ with $\beta-\alpha\in~]0, 2\pi]$ the function $z\longmapsto f(\Re z, \Im z, {\rm arg}_\alpha z)$ is analytic on $S_{\alpha,\beta}.$
\end{Definition}
As an example, the function $L$ defined by
$$ L(r\cos \theta, r\sin \theta, \theta):= \log r+i \theta \quad \quad (r>0,~ \theta \in \mathbb{R}),$$
is analytic on $\mathbb{E}$ and coincides on $\mathbb{R}^+$ with the logarithm of real analysis.

Given $\alpha,$ the largest admissible $\beta$ in the previous definition is $\alpha + 2\pi.$
For $k\in \mathbb{Z},$ all the sectors $S_{\alpha+2k\pi,\alpha +2(k+1)\pi}$ coincide with the complex plane slit along the ray $z=r e^{i \alpha} (r>0)$ and the analytic functions induced by $f$ on these sectors are the analytic branches of $f.$
In this setting, the Paley--Wiener theorem for Mellin band-limited functions in the previous subsection may be restated as follows.
\begin{Theorem}[Paley--Wiener]\label{pw1}
A function $\varphi \in X^2_c$ belongs to the Paley--Wiener space $B^2_{c,T}$ if and only if there exists a function $f: \mathbb{E}\rightarrow \mathbb{C}$ with the following properties:
\begin{enumerate}
\item[(i)] $f$ is analytic on $\mathbb{E};$
\item[(ii)] $f(r,0,0)= \varphi(r)$ \quad $(r>0);$
\item[(iii)] there exists a constant $C_f$ such that
$$ |f(r\cos \theta,r\sin \theta,\theta)| \leq C_f r^{-c} e^{T|\theta|} \quad (\theta \in \mathbb{R});$$
\item[(iv)]$ \lim_{r\rightarrow 0} r^c f(r\cos \theta,r\sin \theta, \theta) = \lim_{r\rightarrow \infty} r^c f(r\cos \theta,r\sin \theta, \theta) =0$
uniformly with respect to $\theta$ on all compact subintervals of $\mathbb{R}.$
\end{enumerate}
\end{Theorem}
The proof is essentially the same as that of Theorem 3 in \cite{BBMS}, using the restrictions of the function $f$ to the sectors
$S_{2k\pi, 2(k+1)\pi}$ as the analytic branches of $f.$


\section{Analytic functions over the polar plane}

In our previous approach with the helicoidal surface we have functions of three variables $x_1, x_2, x_3$ but we use them only as functions of the 
two variables $r$ and $\theta.$ Since there exists a bijection between the helicoidal surface and the right half-plane understood as the set of 
all points $(r, \theta)$ with $r > 0$ and $\theta \in \mathbb{R},$ one may think of considering functions defined on the right half-plane. However, 
these functions will no longer be analytic in the classical sense. They are differentiable and satisfy the Cauchy--Riemann equations transformed 
into polar coordinates. This approach amounts to taking an analytic function, writing its variable in polar coordinates $z = r e^{i \theta}$ and 
treating $(r, \theta)$ as if they were Cartesian coordinates.

Let $\mathbb{H}:= \{(r,\theta) \in \mathbb{R}^+ \times \mathbb{R}\}$ be the right half-plane and let ${\cal D}$ be a domain in $\mathbb{H}.$

\begin{Definition}\footnote{As far as the authors are aware, this modified notion of analyticity arising by treating the polar coordinates as cartesian coordinates, has as yet not been presented.
  Our definition leads naturally to  the classical Cauchy-Riemann
equations when written in their  polar form, often treated
in the literature. Although other mathematicians may have come across this concept too, it seems that it has not been used for practical purposes so far. In Mellin analysis it turns out to be very helpful for an efficient approach, independent of Fourier analysis. In particular it leads to a precise and simple analysis for functions defined over the Riemann surface of the complex logarithm, via the helicoidal surface.}
We say that $f:{\cal D}\rightarrow \mathbb{C}$ is {\rm polar-analytic} on ${\cal D}$ if for any $(r_0, \theta_0) \in {\cal D}$ the limit
$$\lim_{(r,\theta) \rightarrow (r_0, \theta_0)}\frac{f(r, \theta) - f(r_0, \theta_0)}{re^{i\theta} - r_0e^{i\theta_0}} =: (D_{{\rm pol}}f)(r_0, \theta_0)$$
exists and is the same howsoever $(r, \theta)$ approaches $(r_0, \theta_0)$ within ${\cal D}.$
\vskip0,3cm
For a polar-analytic function $f$ we define the polar Mellin derivative as 
$$\Theta_cf(r,\theta) := re^{i \theta}(D_{{\rm pol}}f)(r, \theta) + cf(r,\theta).$$
\end{Definition}
\vskip0,3cm
\newtheorem{Remark}{Remark}
\begin{Remark}
{\rm It can be verified that $f = u + iv$ with $u,v: {\cal D}\rightarrow \mathbb{R}$ is polar-analytic on ${\cal D}$ if and only if $u$ and $v$ have continuous partial derivatives on ${\cal D}$ that satisfy the differential equations}
\begin{eqnarray}\label{CRE}
\frac{\partial u}{\partial \theta} = - r \frac{\partial v}{\partial r}\,,\quad 
\frac{\partial v}{\partial \theta} = r \frac{\partial u}{\partial r}\,.
\end{eqnarray}
\end{Remark}
Note that these equations coincide with the Cauchy-Riemann equations of an analytic function $g$ defined by $g(z) := u(r,\theta) + i v(r, \theta)$ for $z= r e^{i\theta}.$ For the derivative $D_{{\rm pol}}$ we easily find that 
$$(D_{{\rm pol}}f)(r, \theta) = e^{-i\theta}\bigg[\frac{\partial}{\partial r}u(r, \theta) + i \frac{\partial}{\partial r}v (r, \theta) \bigg] = 
\frac{e^{-i\theta}}{r}\bigg[\frac{\partial}{\partial \theta}v (r, \theta) - i \frac{\partial}{\partial \theta}u (r, \theta) \bigg].$$
Also note that $D_{{\rm pol}}$ is the ordinary differentiation on $\mathbb{R}^+.$ More precisely, if $\varphi (\cdot) := f(\cdot, 0)$ then $(D_{{\rm pol}}f)(r,0) =
\varphi'(r).$

Moreover, for $\theta = 0$ we obtain the known formula for $\varphi$ $$\Theta_c\varphi(r) = r\varphi'(r) + c\varphi (r).$$

When $g$ is an entire function, then $f: (r, \theta) \mapsto g(re^{i\theta})$ defines a function $f$ on $\mathbb{H}$ that is polar-analytic 
and $2\pi$-periodic with respect to $\theta.$ The converse is also true. However, there exist polar-analytic functions on $\mathbb{H}$ 
that are not $2\pi$-periodic with respect to $\theta.$ A simple example is the function $L(r, \theta):= \log r + i\theta,$ which is easily seen 
to satisfy the differential equations (\ref{CRE}).

 A connection with the analytic functions in classical sense can be established by a suitable substitution (see the proof of Theorem \ref{asympt} below). Using this, we find that every polar-analytic function $f$ has a series expansion of type
$$f(r,\theta) = \sum_{n=o}^\infty a_n (\log r + i\theta)^n,$$
convergent everywhere on $\mathbb{H}.$ 

Now, for $\alpha, \beta \in \mathbb{R}$ with $\alpha < \beta,$ we consider the set
$$\mathbb{H}_{\alpha, \beta}:= \{(r, \theta) \in \mathbb{R}^+ \times \mathbb{R}: \theta \in ]\alpha, \beta[\}$$
and call it a {\it strip} of $\mathbb{H}.$

If $f : \mathbb{H} \rightarrow \mathbb{C}$ is polar-analytic but not $2\pi$-periodic with respect to $\theta,$ then we can associate with $f$ 
a function $g$ that is analytic on the Riemann surface $S_{{\rm log}}$ of the logarithm. The restriction of $f$ to a strip 
$\mathbb{H}_{\alpha + 2k\pi, \alpha + 2(k+1)\pi},$ where $k \in \mathbb{Z},$ defines an analytic function $g_k$ in the slit complex plane 
$\mathbb{C}\setminus \{re^{i\alpha}: r>0\}$ by setting $g_k(re^{i\theta}):= f(r,\theta).$ The functions $g_k$ for $k \in \mathbb{Z}$ are 
the analytic branches of $g.$
\vskip0,4cm
We now study line integrals for polar-analytic functions. Here a piecewise continuously differentiable curve will be called a {\it regular curve}. The following proposition will be useful in what follows.
\newtheorem{Proposition}{Proposition}
\begin{Proposition}\label{lineintegral}
Let $f$ be a polar-analytic function on $\mathbb{H}$ and let $(r_1, \theta_1)$ and $(r_2, \theta_2)$ be any two points in $\mathbb{H}.$ Then the line integral
\begin{eqnarray}\label{linint}
\int_\gamma f(r,\theta) e^{i\theta}(dr + ird\theta)
\end{eqnarray}
has the same value for each regular curve $\gamma$ in $\mathbb{H}$ that starts at $(r_1,\theta_1)$ and ends $(r_2,\theta_2).$ In particular, the integral vanishes for closed regular curves.
\end{Proposition}
{\it Proof}. Recalling equations (\ref{CRE}), we easily verify that
$$\frac{\partial}{\partial \theta}\bigg[f(r, \theta)e^{i\theta}\bigg] = 
\frac{\partial}{\partial r}\bigg[f(r, \theta) ir e^{i\theta}\bigg].$$
By a theorem of Schwartz, this implies that the integrand in (\ref{linint}) is an exact differential on $\mathbb{H}$, that is, there exists a function $F :\mathbb{H}\rightarrow \mathbb{C}$ such that
$$ \frac{\partial F}{\partial r}(r,\theta) = f(r, \theta)e^{i\theta}\quad \mbox{and} \quad 
\frac{\partial F}{\partial \theta}(r, \theta) = f(r,\theta) ir e^{i\theta},$$
and so the integral in (\ref{linint}) is equal to $F(r_2,\theta_2) - F(r_1, \theta_1).$ \hfill $\Box$
\vskip0,3cm

\section{The Mellin--Bernstein classes}

In this section we introduce Mellin--Bernstein spaces in a somewhat different way closer to the classical one of Fourier analysis.
\begin{Definition}\label{MBS}
For $c \in \mathbb{R},$ $T>0$ and $p \in [1, +\infty[$ the {\rm Mellin--Bernstein space} $\mathscr{B}^p_{c,T}$ comprises all functions $f: \mathbb{H}\rightarrow \mathbb{C}$ with the following properties:
\begin{enumerate}
\item[(i)] $f$ is polar-analytic on $\mathbb{H};$
\item[(ii)] $f(\cdot, 0) \in X^p_{c};$
\item[(iii)] there exists a positive constant $C_f$ such that
$$|f(r,\theta)| \leq C_fr^{-c}e^{T|\theta|}\qquad ((r, \theta) \in \mathbb{H}).$$
\end{enumerate}
\end{Definition}
The above definition has two useful consequences, which are described in the following
\begin{Theorem}\label{asympt}
Let $f \in \mathscr{B}^p_{c,T}$ with $p \in [1, +\infty].$ Then the following statements hold:
\begin{enumerate}
\item[(i)] $f(\cdot, \theta) \in X^p_c$ for all $\theta \in \mathbb{R}$ and $\|f(\cdot, \theta)\|_{X^p_c} \leq e^{T|\theta|}\|f(\cdot, 0)\|_{X^p_c};$
\item[(ii)] $\lim_{r \rightarrow 0}r^cf(r,\theta) = \lim_{r \rightarrow +\infty}r^cf(r,\theta) = 0$ uniformly with respect to $\theta$ on all compact subinterval of $\mathbb{R}.$
\end{enumerate}
\end{Theorem}
{\it Proof}. The function $g: x+iy \rightarrow f(e^x,y)$ is defined on $\mathbb{C}.$ Writing $f = u + iv$ with real-valued functions $u$ and $v,$ 
we know that the differential equations (\ref{CRE}) hold. Therefore $g$ satisfies the Cauchy--Riemann equations on $\mathbb{C}$, and so $g$ is 
an entire function. Consequently
$$F: x+iy \longmapsto e^{c(x+iy)}f(e^x,y)$$
is also an entire function. Property (iii) of Definition \ref{MBS} may now be rewritten in terms of $F$ as
\begin{eqnarray}\label{expo2}
|F(x+iy)| \leq C_fe^{T|y|},
\end{eqnarray}
which shows that $F$ is an entire function of exponential type $T.$ Moreover, property (ii) of Definition 4 shows by a substitution of variables that
$$\int_{-\infty}^{+\infty} |F(x)|^pdx = \|f(\cdot, 0)\|_{X^p_c}^p$$
exists. Now using \cite[Theorem~6.7.1]{BO}, we have
\begin{eqnarray}\label{estim}
\int_{-\infty}^{+\infty} |F(x+iy)|^pdx \leq e^{pT|y|}\int_{-\infty}^{+\infty} |F(x)|^pdx
\end{eqnarray}
for all $y \in \mathbb{R}$ and $F(x) \rightarrow 0$ as $x \rightarrow \pm \infty.$ Therefore assertion (i) follows from (\ref{estim}) by the 
substitution $e^x \mapsto r$ and $y\mapsto \theta.$

Next we note that in view of (\ref{estim}), the convergence $F(x) \rightarrow 0$ extends to $F(x+iy) \rightarrow 0$ as  $x \rightarrow \pm \infty$ 
pointwise for each $y \in \mathbb{R}.$ In connection with (\ref{expo2}), a result in \cite[Theorem 1.4.9]{BO} guarantees that the convergence 
is even uniform with respect to $y$ on compact subintervals of $\mathbb{R};$ also see \cite[p.~170]{TIT}. Writing this statement in terms of $f,$ we find that assertion (ii) is 
also true. \hfill $\Box$

\section{The Paley--Wiener theorem for $\mathscr{B}^2_{c,T}$}

We are ready to state and prove an alternative version of the Paley--Wiener theorem, in terms of the Mellin--Bernstein space $\mathscr{B}^2_{c,T}.$
\begin{Theorem}[Paley-Wiener]\label{pw2}
A function $\varphi \in X^2_c$ belongs to the Paley--Wiener space $B^2_{c,T}$ if and only if there exists a function $f\in \mathscr{B}^2_{c,T}$ such that 
$f(\cdot, 0) = \varphi(\cdot).$
\end{Theorem}
{\it Proof}. 
First, suppose that $\varphi \in B^2_{c,T}.$ Then, by the inversion formula,  we have
$$ \varphi (r) = \frac{1}{2\pi}\int_{-T}^T [\varphi]^\wedge_{M^2_c}(c+it) r^{-c-it}dt\qquad (r>0).$$
We now extend $\varphi$ to the $(r,\theta)$-plane with $r>0$ and $\theta \in \mathbb{R},$ by replacing $r$ with $re^{i\theta}$ on the right-hand side. 
Denoting  this extension by $f(r, \theta)$, we have
\begin{eqnarray}\label{rep}
f(r,\theta) := \frac{1}{2\pi}\int_{-T}^T [\varphi]^\wedge_{M^2_c}(c+it) (re^{i\theta})^{-c-it}dt.
\end{eqnarray}
This imples that 
$$|f(r,\theta)| \leq \frac{r^{-c}e^{T|\theta|}}{2\pi}  \int_{-T}^T |[\varphi]^\wedge_{M^2_c}(c+it) |dt.$$
Hence (iii) of Definition \ref{MBS} holds with 
$$C_f:= \frac{1}{2\pi}  \int_{-T}^T |[\varphi]^\wedge_{M^2_c}(c+it) |dt =
\frac{1}{2\pi}  \int_{-\infty}^{+\infty}|[\varphi]^\wedge_{M^2_c}(c+it) |dt.$$
Concerning polar analyticity of $f,$ for a fixed point $(r_0, \theta_0) \in \mathbb{H},$ we consider the difference quotient
$$\frac{f(r, \theta) - f(r_0, \theta_0)}{re^{i\theta} - r_0e^{i\theta_0}} =  \frac{1}{2\pi}  \int_{-T}^T [\varphi]^\wedge_{M^2_c}(c+it) 
\frac{(re^{i\theta})^{-c-it} - (r_0e^{i\theta_0})^{-c-it}}{re^{i\theta} - r_0e^{i\theta_0}}dt.$$
The limit $(r, \theta) \rightarrow (r_0, \theta_0)$ carried out inside the integral leads to an ordinary differentiation of $z^{-c-it}$ with respect 
to $z$ at the point $z_0:= r_0e^{i\theta_0}.$ We have to justify that the limit and the integration can be interchanged. In order to do that, 
let us take a closed rectangle $Q$ centered at the point $(r_0, \theta_0)$ such that $Q \subset \mathbb{H}.$ We note that the function $h:Q\times [-T,T] \rightarrow \mathbb{C}$, 
defined by
$$ h(r,\theta, t) := \frac{(re^{i\theta})^{-c-it} - (r_0e^{i\theta_0})^{-c-it}}{re^{i\theta} - r_0e^{i\theta_0}}~~\mbox{for}~(r,\theta) \neq (r_0, \theta_0),$$
and
$$ h(r_0,\theta_0,t):= -(c+it) (r_0e^{i\theta_0})^{-c-1-it},$$
is continuous and therefore its absolute value is bounded on $Q\times [-T,T].$ Hence the desired interchange is guaranteed by Lebesgue's theorem of
dominated convergence.

Now we prove the reverse implication.  Let $f\in \mathscr{B}^2_{c,T}$. For any fixed $t \in \mathbb{R},$  define the function
$$g(r, \theta) := e^{(c-1+it)(\log r +i\theta)}f(r, \theta),$$
which is polar-analytic in $\mathbb{H}.$ 

Let us first assume that $t>T.$ For $R>1$, we define $n:=\lfloor \log R\rfloor$, where for any real number $a$ the symbol $\lfloor a\rfloor$ denotes the
largest integer not exceeding $a$.  We now consider the regular curve $\gamma$ which is the boundary of the rectangle 
in $\mathbb{H}$ defined by $[1/R, R] \times [0, 2n\pi].$ Then, using Proposition \ref{lineintegral}, one has
$$\int_\gamma g(r,\theta) e^{i\theta}(dr + ird \theta) = 0.$$
This easily implies that 
\begin{eqnarray*}
\lefteqn{\int_{1/R}^R r^{c+it -1}\varphi(r)dr = \int_{1/R}^R e^{(c-1+it)(\log r + 2n\pi i)}f(r, 2n\pi)dr } \\ 
& & +\frac{i}{R}\int_0^{2n \pi} e^{(c-1+it)(-\log R +i\theta) + i\theta}f\left(\frac{1}{R}, \theta\right) d\theta -iR\int_0^{2n \pi}e^{(c-1+it)(\log R +i\theta) + i\theta}f(R, \theta)d\theta\\
&=:& I_1 + I_2 -I_3,
\end{eqnarray*}
where $\varphi$ is the restriction of $f$ to $\R^+$.

We now estimate the integrals on the right-hand side.
Using (iii) of Definition~\ref{MBS}, we find that
\begin{eqnarray*}
|I_1| &\leq& \int_{1/R}^R e^{(c-1)\log r}e^{-2n\pi t}|f(r, 2\pi n)|dr \\ &\leq& 
C_f\int_{1/R}^R e^{-(t-T)2n\pi}\,\frac{dr}{r} = C_fe^{-(t-T)2n\pi} \log R^2.
\end{eqnarray*}
Thus, according to our choice of $n$, we obtain $|I_1|\rightarrow 0$ as $R\rightarrow +\infty.$

As to $I_2,$ we easily see that
$$|I_2| \leq \int_0^{2n\pi} \frac{1}{R^c}\bigg|f\left(\frac{1}{R}, \theta\right)\bigg| e^{-t\theta}d\theta \leq
\int_0^{+\infty} \frac{1}{R^c}\bigg|f\left(\frac{1}{R}, \theta\right)\bigg| e^{-t\theta}d\theta.
$$
Now, since  the integrand is dominated by $C_f e^{-(t-T)\theta}$, which is integrable as a function of $\theta$, by (ii) of Theorem \ref{asympt} 
and the theorem of dominated convergence, we obtain $|I_2| \rightarrow 0$ as $R\rightarrow +\infty.$ The same result is obtained for $|I_3|$. 
Altogether we have
$$\lim_{R \rightarrow +\infty}\left|\int_{1/R}^R r^{c-1+it}\varphi(r)dr\right| = 0.$$
This implies that $[\varphi]^\wedge_{M^2_c}(c+it) = 0.$ 

When $t<-T$, we choose $n:= -\lfloor\log R\rfloor$ and proceed analogously using the boundary of the rectangle $[1/R, R]\times [-2n\pi, 0].$ 
Hence $\varphi(\cdot):= f(\cdot, 0)$ belongs to $B^2_{c,T}.$
\hfill $\Box$
\vskip0,4cm

\section{A Hardy space in Mellin analysis}

In this section we apply our approach in order to define in a simple way a Hardy space in the Mellin frame. For $c \in \mathbb{R}$ and  
$p \in [1,+\infty[$,  we recall that the norm 
in $X^p_c$ is defined by
$$\|\varphi\|_{X^p_c}= \bigg(\int_0^{+\infty} |\varphi(r)|^p r^{cp-1}dr\bigg)^{1/p}.$$
For a strip $\mathbb{H}_{-a,a}$ with $a>0$, we simply write $\mathbb{H}_a.$ 
\begin{Definition}\label{Hardy}
Let $a, c, p \in \mathbb{R}$ with $a>0$ and $p\geq 1.$ The {\rm Mellin--Hardy space} $H^p_c(\mathbb{H}_a)$ comprises all functions $f : \mathbb{H}_a \rightarrow \mathbb{C}$ that satisfy the following 
conditions:
\begin{enumerate}
\item[(i)] $f$ is polar-analytic on $\mathbb{H}_a;$
\item[(ii)] $f(\cdot, \theta) \in X^p_c$ for each $\theta \in {]}-a,a{[};$
\item[(iii)] there holds
$$\|f\|_{H^p_c(\mathbb{H}_a)} := \sup_{0<\theta<a}\bigg(\frac{\|f(\cdot, \theta)\|^p_{X^p_c} + \|f(\cdot, -\theta)\|^p_{X^p_c}}{2}\bigg)^{1/p} < +\infty.$$
\end{enumerate}
\end{Definition}
When $a \in {]}0,\pi]$ we can associate with each function $f \in \mathbb{H}_a$ a function $g$ analytic on the sector 
${\cal S}_a:= \{z \in \mathbb{C}: |\arg z| < a\}$ by defining 
$g(re^{i\theta}):= f(r,\theta).$ The collection of all such functions constitutes a Hardy-type space $H^p_c({\cal S}_a),$ which may be identified 
with $H^p_c(\mathbb{H}_a).$

In \cite{BSS2} the authors considered a Hardy space $H^p(S_a)$ of functions analytic on the strip $S_a:=\{z \in \mathbb{C}: |\Im z| < a\}.$ The Hardy space $H^p_c(\mathbb{H}_a)$ has been designed in such way that if $g \in H^p(S_a)$ and 
$$f(r,\theta):= r^{-c}e^{-ic\theta}g(\log r + i\theta),$$
then $f \in H^p_c(\mathbb{H}_a)$ and conversely, if $f \in H^p_c(\mathbb{H}_a)$ and 
$$g(x+iy):= e^{c(x+iy)}f(e^x,y),$$
then $g \in H^p(S_a).$

Using this correspondence, one can deduce the following propositions from the analogous ones proved in \cite{BSS2}. 
\begin{Proposition}\label{p1}
Let $f \in H^p_c(\mathbb{H}_a)$ and let $a_1 \in {]}0,a{[}.$ Then
$$|f(r,\theta)| \leq r^{-c}\bigg(\frac{4}{\pi(a-a_1)}\bigg)^{1/p}\|f\|_{H^p_c(\mathbb{H}_a)}$$
for $(r,\theta) \in \mathbb{H}_{a_1}.$
\end{Proposition}
The next proposition is a Nikol'ski-type inequality for the Hardy space $H^p_c(\mathbb{H}_a).$
\begin{Proposition}\label{p2}
Let $f \in H^p_c(\mathbb{H}_a)$ and let $(t_n)_{n \in \mathbb{Z}}$ be a sequence on $\mathbb{R}^+$ such that $t_{n+1}/t_n > e^{2\delta}$ for all 
$n \in \mathbb{Z}$ with 
$\delta \in {]}0,a{[}.$ Then
$$\bigg(\sum_{n \in \mathbb{Z}} t_n^{cp}|f(t_n,0)|^p\bigg)^{1/p} \leq \bigg(\frac{2}{\pi \delta}\bigg)^{1/p} \|f\|_{H^p_c(\mathbb{H}_a)}.$$
\end{Proposition}
\begin{Proposition}\label{p3}
Let $f \in H^p_c(\mathbb{H}_a)$ and let $a_1 \in {]}0,a{[}.$ Then
$$\lim_{r \rightarrow 0^+}r^cf(r,\theta) =\lim_{r \rightarrow +\infty}r^cf(r,\theta) =0$$
uniformly for $\theta \in [-a_1, a_1].$
\end{Proposition}
The following theorem will be needed for estimating the distance of a function $f\in H^p_c(\mathbb{H}_a)$ from the Mellin--Bernstein space 
$\mathscr{B}^p_{c,\sigma}$ (see the next section). Here we will restrict ourselves to $p\in \{1,2\}.$
\begin{Theorem}\label{mellintransf}
Let $p \in \{1,2\}$ and $f\in H^p_c(\mathbb{H}_a).$ Then for $\alpha \in {]}0,a{[}$, we have
$$|[f(\cdot,0)]^\wedge_{M^p_c}(c+it)| = e^{-\alpha |t|}|[f(\cdot, \varepsilon \alpha)|^\wedge_{M^p_c}(c+it)| \qquad (a.e.~ for~ t \in \mathbb{R}),$$
where $\varepsilon = 1$ for $t>0$ and $\varepsilon = -1$ for $t<0.$
\end{Theorem}
{\it Proof}. For $f \in H^p_c(\mathbb{H}_a)$ let us consider the function $g$ defined by
$$g(r,\theta):= e^{(c-1+it)(\log r + i\theta)}f(r, \theta),$$
which is polar-analytic on $\mathbb{H}_a.$ For  sake of simple notation, we do not indicate  the dependence of $g$ on $c$ and $t$ explicitely.

First we assume that $t>0.$ For $\rho >1,$ let $\gamma$ be the positively oriented rectangular curve in $\mathbb{H}$ with vertices at 
$(1/\rho, 0), (\rho, 0), (1/\rho, \alpha)$ and $(\rho, \alpha).$ By Proposition \ref{lineintegral}, we have
$$\int_\gamma g(r,\theta) e^{i\theta}(dr + ird\theta) = 0.$$
This equation may be rewritten in terms of ordinary integrals as
$$\int_{1/\rho}^\rho g(r,0)dr = \int_{1/\rho}^\rho g(r,\alpha) e^{i\alpha}dr + I_+\bigg(\frac{1}{\rho}, t\bigg) - I_+(\rho, t),$$
where
$$I_+(r,t) = \int_0^\alpha g(r,\theta)i r e^{i\theta}d\theta = i r^{c+it}\int_0^\alpha e^{-(t-ic)\theta}f(r, \theta)d\theta,$$
and so 
$$|I_+(r,t)| \leq r^c\int_0^\alpha e^{-t \theta}|f(r,\theta)| d\theta.$$
Now, for $p=2$ the integrals
$$\int_{1/\rho}^\rho g(r,0)dr \quad\hbox{and} \quad \int_{1/\rho}^\rho g(r,\alpha) e^{i\alpha}dr$$
converge in the $L^2$-sense to the Mellin transforms 
$$[f(\cdot, 0)]^\wedge_{M^p_c}(c+it)  \quad \hbox{and} \quad 
e^{-\alpha t} e^{i \alpha c}[f(\cdot, \alpha)]^\wedge_{M^p_c}(c+it),$$ 
respectively, while for $p=1$ they converge uniformly.

 For $t<0,$ we use the rectangular curve $\gamma$ with vertices at $(1/\rho, 0), (\rho, 0), (1/\rho, -\alpha)$ and $(\rho, -\alpha),$ and 
proceed analogously. In any case, we obtain for all $t \in \mathbb{R},$
\begin{eqnarray}\label{gamma}
\int_{1/\rho}^{\rho} r^{c-1+it}f(r,0)dr = e^{-\alpha |t|} e^{i\varepsilon \alpha c}\int_{1/\rho}^\rho r^{c-1+it}f(r, \varepsilon \alpha) dr + 
I\bigg(\frac{1}{\rho}, t\bigg) - I(\rho, t),
\end{eqnarray}
where
$$|I(r,t)| \leq r^c\int_{-\alpha}^\alpha e^{-|\theta t|}|f(r, \theta)| d\theta.$$
Using Proposition \ref{p3}, we can show that $\|I(r, \cdot)\|_{L^{p'}(\mathbb{R})}$ exists for $p' \in \{2, \infty\}$ and approaches $0$ as 
$r \rightarrow 0$ or $r \rightarrow +\infty.$ Hence (\ref{gamma}) implies the assertion. \hfill $\Box$

\section{Estimates of the distance from Mellin--Bernstein spaces}

In \cite{BBMS2} we introduced a notion of distance in terms of the Mellin transform. For $c \in \mathbb{R}$ and $q \in [1, +\infty],$ let $G^q_c$ 
be the linear space of all functions $f:\mathbb{R}^+\rightarrow \mathbb{C}$ that have the representation
$$f(x) = \frac{1}{2\pi}\int_{-\infty}^\infty \psi(v) x^{-c-iv}dv \qquad (x>0),$$
where $\psi \in L^1(\mathbb{R})\cap L^q(\mathbb{R}).$ We endowed this space with the norm
$$[\!\![f]\!\!]_q : = \|\psi\|_{L^q(\mathbb{R})} = \bigg(\int_{\mathbb{R}}|\psi(v)|^qdv\bigg)^{1/q}$$
and the corresponding metric is defined by
$$\mbox{dist}_q(f,g) := [\!\![f-g]\!\!]_q, \qquad f,g \in G^q_c.$$
Note that if $f \in \mathscr{B}^p_{c,\sigma},$ for $p \in [1,2],$ then $f(\cdot, 0) \in G^q_c$ for any $q \geq 1$, since its 
Mellin transform has compact support, and therefore, from the uniqueness theorem of Mellin transforms, one can take $\psi (v)= [f(\cdot,0)]^\wedge_{M^p_c}(c+iv).$

In polar form, we can define the class $\widetilde{G}^q_c$ as the space of all polar-analytic functions on $\mathbb{H}$ such that $f(\cdot, 0) \in G^q_c.$

Another important subspace is the {\it Mellin inversion class}, denoted by ${\cal M}^p_c,$
with $p \in [1,2],$ consisting of all functions $f \in X^p_c \cap C(\mathbb{R}^+)$ such that $[f]^\wedge_{M^p_c} \in L^1(\{c\} \times i\mathbb{R}).$ 
This class is contained in $G^q_c$ for $q \in [1, p']$ with $p'$ being the conjugate exponent of $p.$

For $p \in \{1,2\}$, the Mellin--Hardy space defined above contains the space $\mathscr{B}^p_{c,\sigma}.$ Moreover, using Proposition \ref{p1} and 
Theorem \ref{mellintransf}, we see that the restriction to the positive real line of a function in the Hardy space is an element of ${\cal M}^p_c.$

For a subspace $A \subset G^q_c$ we define
$$\mbox{dist}_q(f, A) := \inf_{g \in A}[\!\![f-g]\!\!]_q.$$
In \cite{BBMS2} we obtained a representation formula for the distance of a function $f \in G^q_c$ from the Paley--Wiener space 
$B^p_{c,\sigma}$  in the form
$$\mbox{dist}_q(f , B^p_{c,\sigma}) = \bigg(\int_{|v| \geq \sigma} |\psi(v)|^q dv\bigg)^{1/q} \qquad (1 \leq q < \infty)$$
and if $\psi$ is also continuous, then
$$\mbox{dist}_\infty(f, B^p_{c,\sigma}) = \sup_{|v|\geq \sigma}|\psi(v)|.$$
Moreover, in the same paper, we estimated the distance in case of Mellin inversion classes, Lipschitz spaces and Mellin-Sobolev spaces.

Now, we define the distance $\mbox{dist}_q(f, \mathscr{B}^p_{c, \sigma})$ of a function $f \in H^p_c(\mathbb{H}_a),$  with $a >0,$ from the space $\mathscr{B}^p_{c, \sigma},$ by considering the restrictions to the positive real line of the function $f$ and the Mellin-Bernstein spaces $\mathscr{B}^p_{c, \sigma}.$

The following result gives estimates of the Mellin distance $\mbox{dist}_q(f, \mathscr{B}^p_{c, \sigma})$ for a function 
$f \in H^p_c(\mathbb{H}_a)$ with $a>0.$
\begin{Theorem}\label{distance1}
Let $f \in H^p_c(\mathbb{H}_a)$ with $p\in \{1,2\}.$ Then for $q \in [1, \infty],$ we have
\begin{eqnarray*}
\mbox{\rm dist}_q(f, \mathscr{B}^1_{c, \sigma}) & \leq &2\|f\|_{H^1_c(\mathbb{H}_a)} \bigg(\frac{2}{aq}\bigg)^{1/q} e^{-a\sigma}\\
\mbox{\rm dist}_2(f, \mathscr{B}^2_{c, \sigma})  &\leq &2\sqrt{\pi}\|f\|_{H^2_c(\mathbb{H}_a)} e^{-a\sigma}.
\end{eqnarray*}
Furthermore,  for $q \in [1,2{[}$,
$$\mbox{\rm dist}_q(f, \mathscr{B}^2_{c, \sigma}) \leq 2\sqrt{\pi}\bigg(\frac{2-q}{qa}\bigg)^{(1/q) - (1/2)}\|f\|_{H^2_c(\mathbb{H}_a)} e^{-a\sigma}.$$
\end{Theorem}
{\it Proof}. For $p=1,$ let us consider the function $f(\cdot, \alpha \varepsilon)$ as in Theorem \ref{mellintransf}. It is easy to see that
$$|[f(\cdot,\alpha \varepsilon)]^\wedge_{M_c}(c+it)| \leq \|f(\cdot,\alpha\varepsilon)\|_{X_c} \leq 2\|f\|_{H^1_c(\mathbb{H}_a)}.$$
Then, by Theorem \ref{mellintransf},
\begin{eqnarray*}
\lefteqn{\mbox{dist}_q(f, \mathscr{B}^1_{c, \sigma}) = \bigg(\int_{|t| \geq \sigma}|[f(\cdot, 0)]^\wedge_{M_c}(c+it)|^qdt\bigg)^{1/q}}\\
&=& \bigg(\int_{-\infty}^{-\sigma} e^{-q\alpha |t|}\abs{[f(\cdot, -\alpha)]^\wedge_{M_c}(c+it)}^q dt
     + \int_\sigma^\infty e^{-q\alpha |t|}\abs{[f(\cdot, \alpha)]^\wedge_{M_c}(c+it)}^q dt \bigg)^{1/q}\\
&\le& 2 \li\|f\ri\|_{H_c^1(\mathbb{H}_a)} \li(2\int_\sigma^\infty e^{-q\alpha t} dt\ri)^{1/q} \,=\, 
      2 \li\|f\ri\|_{H_c^1(\mathbb{H}_a)} \li(\frac{2}{\alpha q}\ri)^{1/q} e^{-\alpha \sigma}
\end{eqnarray*}
for every $\alpha \in {]}0,a{[}.$ Now the assertion follows for $\alpha \rightarrow a^-.$

Next, let $p=2.$ For $q=2$ we may employ Theorem~\ref{mellintransf} in conjunction with Lemma~\ref{Plancherel} to conclude that
\begin{eqnarray*}
\mbox{dist}_2(f, \mathscr{B}^2_{c, \sigma}) &=&\bigg(\int_{|t| \geq \sigma}\abs{[f(\cdot, 0)]^\wedge_{M^2_c}(c+it)}^2dt\bigg)^{1/2}\\
  &=&\bigg(\int_{|t| \geq \sigma}e^{-2\alpha\abs{t}}\abs{[f(\cdot, \eps\alpha)]^\wedge_{M^2_c}(c+it)}^2dt\bigg)^{1/2}\\
  &\le&e^{-\alpha\sigma}\bigg(\int_{-\infty}^{-\sigma}\abs{[f(\cdot, -\alpha)]^\wedge_{M^2_c}(c+it)}^2dt \\
  && \qquad\quad +\int_\sigma^\infty \abs{[f(\cdot, \alpha)]^\wedge_{M^2_c}(c+it)}^2dt \bigg)^{1/2}\\
  &\le& e^{-\alpha\sigma} \sqrt{2\pi} \li(\li\|f( \cdot, -\alpha)\ri\|_{X_c^2}^2 +\li\|f( \cdot, \alpha)\ri\|_{X_c^2}^2\ri)^{1/2}\\
  &\le& e^{-\alpha\sigma} 2\sqrt{\pi} \li\|f\ri\|_{H_c^2(\mathbb{H}_a)},  
\end{eqnarray*}
and the assertion follows for $\alpha \rightarrow a^-.$

Finally, for $q \in [1,2{[},$ 
we modify the previous proof by using  H\"older's inequality with the conjugate exponents $\mu = 2/(2-q)$ and $\nu = 2/q$.
This leads us to 
\begin{eqnarray*}
\mbox{dist}_q(f, \mathscr{B}^2_{c, \sigma})&=&\bigg(\int_{|t| \geq \sigma}e^{-q\alpha\abs{t}}\abs{[f(\cdot, \eps\alpha)]^\wedge_{M^2_c}(c+it)}^qdt\bigg)^{1/q}\\
  &\le& \bigg(\int_{|t| \geq \sigma} e^{-q\mu \alpha |t|}dt\bigg)^{1/(q\mu)}
 \bigg(\int_{|t| \geq \sigma}\abs{[f(\cdot, \alpha \varepsilon)]^\wedge_{M^2_c}(c+it)}^{q\nu} dt\bigg)^{1/(q\nu)}.
\end{eqnarray*}
The first integral on the right-hand side can be easily calculated. The integrand of the second integral has the exponent $2$. Therefore, we can use
Lemma~\ref{Plancherel} as in the previous paragraph. Thus we obtain
$$\mbox{dist}_q(f, \mathscr{B}^2_{c, \sigma})\,\le\,2\sqrt{\pi}\li\|f\ri\|_{H^2_c(\mathbb{H}_a)}\bigg(\frac{2-q}{q\alpha}\bigg)^{(1/q) - (1/2)}e^{-\alpha \sigma},$$
and the assertion follows again by letting $\alpha \rightarrow a^-.$ \hfill $\Box$

\section{Mellin--Hardy spaces and sampling}
In this section we apply the estimates obtained in the previuos section in order to study the remainder of the exponential sampling formula for 
functions $f$ belonging to a Mellin--Hardy space. We recall here some notions concerning the exponential sampling; see \cite{BJ3}, \cite{BBM1}, \cite{SCH}.

In the following, for $c \in \mathbb{R}$, we denote by $\mbox{lin}_c$  the function
$$\mbox{lin}_c(x) := \frac{x^{-c}}{2\pi i}\frac{x^{\pi i} -x^{-\pi i}}{\log x} = \frac{x^{-c}}{2\pi}\int_{-\pi}^\pi x^{-it}dt \qquad (x>0, \, x\neq 1)$$
with the continuous extension $\mbox{lin}_c(1) = 1.$ Thus 
$$\mbox{lin}_c(x) = x^{-c}\mbox{sinc}(\log x) \qquad (x>0).$$
Here, as usual, the sinc function is defined by 
$$\mbox{sinc}(t):= \frac{\sin (\pi t)}{\pi t}~ \mbox{ \,for }~t\neq 0,\qquad \mbox{sinc}(0) = 1.$$
It is clear that $\mbox{lin}_c \not \in X^1_b$ for any $b\in \R.$ However, it belongs to the space $X^2_c$ and its Mellin transform in 
$X^2_c$-sense is given by
$$[\mbox{\rm lin}_c]^\wedge_{M_c^2}(c+iv) = \chi_{[-\pi, \pi]}(v),$$
where $\chi_A$ denotes the characteristic function of the set $A.$

For a function  $g \in \mathscr{B}^2_{c,\pi T}$ the following exponential sampling formula holds for the restriction of $g$ to the positive real axis 
(see \cite{BJ3}, \cite{BJ2}):
$$g(x) = \sum_{k \in \mathbb{Z}}g(e^{k/T})\mbox{\rm lin}_{c/T}(e^{-k}x^T) \qquad (x>0).$$
As an approximate version in the space $\mathcal{M}^2_c$, we have (see  \cite[Theorem~5.5]{BJ4}):
\begin{Proposition}\label{approsamp}
Let $g \in \mathcal{M}_c^2.$
Then there holds the error estimate
\begin{eqnarray*}
\lefteqn{\bigg|g(x) - \sum_{k=-n}^{n} g(e^{k/T})\mbox{\rm lin}_{c/T}(e^{-k}x^T)\bigg|} \quad\\
&\leq& \frac{x^{-c}}{\pi}\int_{|t| > \pi T}| [g]^\wedge_{M_c^2}(c+it)| dt \qquad(x \in \mathbb{R}^+,~T >0).
\end{eqnarray*}
\end{Proposition}
Note that in general the series in Proposition \ref{approsamp} converges in the principal Cauchy value sense. However, if 
$f \in H^2_c(\mathbb{H}_a),$ then by Proposition \ref{p2} one can deduce that the series for $g(\cdot) := f(\cdot, 0)$ is absolutely and uniformly convergent over any compact subsets of $\mathbb{R}.$ In this case, 
introducing a remainder $(R_{\pi T} g)(x)$ by writing
\begin{eqnarray*}
 g(x) = \sum_{k=-n}^n g\left(e^{k/T}\right) \mbox{lin}_{c/T}\left(e^{-k}x^T\right)
+ \left(R_{\pi T}g\right)(x),
\end{eqnarray*}
we have  by Proposition~\ref{approsamp} in conjunction with Theorem~\ref{pw2}
$$
|(R_{\pi T} g)(x)|\,\leq\, \frac{x^{-c}}{\pi}\,\mbox{dist}_1(g, B_{c,\pi
T}^2) \qquad (x>0),$$
or equivalently,
\begin{eqnarray*}
\|R_{\pi T} g\|_{X_c^\infty}\,\leq\, \frac{1}{\pi}\,\mbox{dist}_1(g, B_{c,\pi T}^2).
\end{eqnarray*}

Thus, as a consequence of Theorem~\ref{distance1}, we obtain:
\newtheorem{Corollary}{Corollary}
\begin{Corollary}
Let $f \in H^2_c(\mathbb{H}_a).$ Then we have
$$\|R_{\pi T} f(\cdot, 0)\|_{X_c^\infty} = {\cal O}(e^{-a\pi T}) \qquad (T \rightarrow +\infty).$$
\end{Corollary}
\vskip0,4cm
Other approximate formulae can be estimated analogously for Mellin--Hardy spaces, as an example, the reproducing kernel formula and a sampling 
formulae for Mellin derivatives; see \cite{BBMS2}.


\section{Acknowledgments}
Carlo Bardaro and Ilaria Mantellini have partially supported by the ``Gruppo Nazionale per l'Analisi Matematica e Applicazioni'' (GNAMPA) of the 
``Istituto Nazionale di Alta Matematica'' (INDAM) as well as by the Department of Mathematics and Computer Sciences of the University of Perugia.


\begin{thebibliography}{99}
\small
\bibitem{AN} N.B. Andersen, On real Paley-Wiener theorems for certain integral transforms, J. Math. Anal. Appl., 288, (2003), 124--135.
\bibitem{AD} N.B. Andersen and M. de Jeu, Real Paley-Wiener theorems and local spectral radius formula, Trans. Amer. Math. Soc., 362,
(2010), 3613--3640.
\bibitem{BBM1} C. Bardaro, P.L. Butzer and I. Mantellini, The exponential sampling theorem of signal analysis and the reproducing kernel formula in the 
Mellin transform setting, Sampl. Theory Signal Image Process., 13(1), (2014), 35--66.
\bibitem{BBM3} C. Bardaro, P.L. Butzer and I. Mantellini, The Mellin-Parseval formula and its interconnections with the exponential sampling theorem of 
optical physics, Integral Transforms and Special Functions, 27(1), (2016), 17--29.
\bibitem{BBMS} C. Bardaro, P.L. Butzer, I. Mantellini and G. Schmeisser, On the Paley-Wiener theorem in the Mellin transform setting, 
J. Approx. Theory, 207, (2016), 60--75.
\bibitem{BBMS2} C. Bardaro, P.L. Butzer, I. Mantellini and G. Schmeisser, Mellin analysis and its basic associated metric. Applications to sampling theory, 
Analysis Math., 42(4), (2016), 297--321.
\bibitem{BO} R.P. Boas, Entire Functions, Academic Press, New York, 1954.
\bibitem{BJ2} P.L. Butzer and S. Jansche,  A direct approach to the Mellin transform, J. Fourier Anal. Appl., 
3, (1997), 325--375.
\bibitem{BP} M. Bertero and E.R. Pike, Exponential sampling method for Laplace and other dilationally invariant transforms I. Singular-system analysis. II. 
Examples in photon correction spectroscopy and Frauenhofer diffraction, Inverse Problems, 7, (1991), 1--20; 21--41.
\bibitem{BJ3} P.L. Butzer and S. Jansche, The exponential sampling theorem of signal analysis, 
Atti Sem. Mat. Fis. Univ. Modena, Suppl. Vol. 46, (1998), 99--122.
\bibitem{BJ4} P.L. Butzer and S. Jansche, A self-contained approach to Mellin transform analysis for square integrable functions, applications, 
Integral Transforms Spec. Funct., 8, (1999), 175--198.
\bibitem{BSS} P.L. Butzer, G. Schmeisser and R.L. Stens, The classical and approximate sampling theorems and their equivalence for entire functions 
of exponential type, J. Approx. Theory, 179, (2014), 94--111.
\bibitem{BSS2}  P.L. Butzer, G. Schmeisser and R.L. Stens, Basic relations valid for the Bernstein spaces $B^2_\sigma$ and their extensions to larger 
functions spaces via a unified distance concept, J. Fourier Anal. Appl., 19, (2013), 333--375.
\bibitem{CBV} R.V. Churchill, J.W. Brown and R.F. Verhey, Complex Variables and Applications (3.~Ed.), McGraw-Hill, New York, 1974.
\bibitem{FS} G.B. Folland and E.M. Stein,  Hardy Spaces on Homogeneous Groups, Mathematical Notes 28, 
Princeton University Press, Princeton, N.J., 1982.
\bibitem{GPS} H-J. Glaeske, A.P. Prudnikov and K.A. Skornik, Operational Calculus and Related Topics, Chapman and Hall, CRC, Boca Raton, FL, 2006.
\bibitem{HIG} J. R. Higgins, Sampling Theory in Fourier and Signal Analysis. Foundations, Oxford Univ. Press, Oxford, 1996.
\bibitem{MA} R.G. Mamedov, The Mellin Transform and Approximation Theory, (in Russian), "Elm", Baku, 1991.
\bibitem{NP} R. Nevanlinna and V. Paatero, Introduction to Complex Analysis, Addison-Wesley Publ. Co., London, 1969.
\bibitem{RS} Q.I. Rahman and G. Schmeisser, $L^p$ inequalities for entire functions of exponential type, Trans. Amer. Math. Soc., 320(1), (1990), 91--103.
\bibitem{RR} M. Rosenblum and J. Rovnyak, Hardy Classes and Operator Theory, Oxford Mathematical Monographs,  
The Clarendon Press, Oxford, 1985.
\bibitem{RU} W. Rudin, Real and Complex Analysis (3.~Ed.), McGraw-Hill, New York, 1986.
\bibitem{RU1} W. Rudin, Functional Analysis (2.~Ed.), McGraw-Hill, New York, 1991.
\bibitem{SCH} G. Schmeisser, Quadrature over a semi-infinite interval and Mellin transform, in: Y.~Lyubarskii (ed.) ``Proceedings of the
1999 International Workshop on Sampling Theory and Applications'', 
(ISBN 82-7151-0991), Norwegian University of 
Science and Technology, Trondheim, 1999, pp.~203--208.
\bibitem{TIT} E.C. Titchmarsh, Introduction to the Theory of Fourier Integrals (2.~Ed.), Clarendon Press, Oxford, 1948.
\bibitem{TU} V.K. Tuan, New type Paley--Wiener theorems for the modified multidimensional Mellin transform, J. Fourier Anal. Appl., 4, (1998), 317--328.
\bibitem{YO} K. Yosida, Functional Analysis (6.~Ed.), Springer-Verlag, Berlin, 1980.


\end{thebibliography}
\end{document}